\documentclass[11pt,a4paper]{book}
\usepackage{jfaa2eplain}

\usepackage{amsmath}
\usepackage[isolatin]{inputenc}
\usepackage{graphicx}
\usepackage{pgf}
\usepackage{cite}

\DeclareMathOperator{\sgn}{sgn}

\DeclareMathOperator*{\argmin}{argmin}
\DeclareMathOperator{\rg}{rg}

\font\msbm=msbm10
\def\mathbb#1{\hbox{\msbm{#1}}}
\newcommand{\RR}{\mathbb{R}}
\newcommand{\NN}{\mathbb{N}}
\newcommand{\HH}{\mathcal{H}}
\newcommand{\lpspace}[1]{\ell^{#1}}
\newcommand{\norm}[2][]{\|{#2}\|_{#1}}
\newcommand{\abs}[1]{{|{#1}|}}
\newcommand{\scp}[3][]{\langle{#2},\, {#3}\rangle_{#1}}
\newcommand{\sign}[1]{\sgn({#1})}
\newcommand{\sett}[1]{\{{#1}\}}
\newcommand{\seq}[1]{\{{#1}\}}
\newcommand{\subgrad}{\partial}
\newcommand{\dd}[1]{\ \mathrm{d}{#1}}
\newcommand{\bigabs}[1]{\bigl|{#1}\bigr|}
\newcommand{\Bigabs}[1]{\Bigl|{#1}\Bigr|}
\newcommand{\placeholder}{\,\cdot\,}
\newcommand{\set}[2]{\{{#1} \ : \ {#2}\}}
\newcommand{\without}{\backslash}
\newcommand{\bigscp}[3][]{\bigl\langle{#2},\, {#3}\bigr\rangle_{#1}}
\newcommand{\Bigset}[2]{\Bigl\{{#1} \ : \ {#2}\Bigr\}}
\newcommand{\kronO}{\mathcal{O}}
\newcommand{\kernel}[1]{\ker({#1})}
\newcommand{\range}[1]{\rg({#1})}

\begin{document}
\author{Kristian Bredies and Dirk A.~Lorenz}

\chapter{Linear convergence of iterative soft-thresholding}

\footnotetext{\textit{Math Subject Classifications.}
                      65J22, 46N10, 49M05.
%
%
}

\footnotetext{\textit{Keywords and Phrases.}
                     Iterative soft-thresholding, inverse problems,
                     sparsity constraints, convergence analysis,
                     generalized gradient projection method.}

\begin{abstract}
  In this article a unified approach to iterative soft-thresholding algorithms for the solution of linear operator equations in infinite dimensional Hilbert spaces is presented.
  We formulate the algorithm in the framework of generalized gradient methods and present a new convergence analysis.
  As main result we show that the algorithm converges with linear rate
  as soon as the underlying operator satisfies the so-called finite
  basis injectivity property or the minimizer possesses a
  so-called strict sparsity pattern.
  Moreover it is shown that the constants can be calculated explicitly
  in special cases (i.e.~for compact
  operators). Furthermore, the techniques also can be used to establish
  linear convergence for related methods such as the 
  iterative thresholding algorithm for joint sparsity and 
  the accelerated gradient projection method.
  
\end{abstract}

\section{Introduction}
\label{sec:intro}

This paper is concerned with the convergence analysis of numerical algorithms
for the solution of linear inverse problems in the
infinite-dimensional setting  with so-called sparsity constraints.
The background for this type of problem is, for example, the attempt
to solve the linear operator equation $Ku = f$ in an
infinite-dimensional Hilbert space which models the connection between
some quantity of interest $u$ and some measurements
$f$. Often, the measurements $f$ contain noise which makes the
direct inversion ill-posed and practically impossible. Thus, instead
of considering the linear equation, a regularized problem is posed for
which the solution is stable with respect to noise.
A common approach is to regularize by minimizing a
Tikhonov functional
\cite{engl1996inverseproblems,daubechies2003iteratethresh,
  resmerita2005regbanspaces}.
 A special class of these regularizations has been
of recent interest, namely of the type
\begin{equation}
  \label{eq:sparse_min_prob}
  \min_{u \in \lpspace{2}} \  \frac{\norm{Ku - f}^2}{2} +
  \sum_{k=1}^\infty \alpha_k \abs{u_k} \ .
\end{equation}
These problems model the fact that the quantity of interest $u$ is
composed of a few elements, i.e.~it is sparse in some given, countable
basis.
To make this precise, let $A:\HH_1\to\HH_2$ be a bounded operator between two Hilbert spaces and let $\{\psi_k\}$ be an orthonormal basis of $\HH_1$.
Denote with $B:\lpspace{2}\to\HH_1$ the synthesis operator
$B(u_k) = \sum_{k}u_k\psi_k$. Then the problem
\[
\min_{u\in\HH_1}\frac{\norm{Au - f}^2}{2} +
  \sum_{k=1}^\infty \alpha_k \abs{\scp{u}{\psi_k}}
\]
can be rephrased as \eqref{eq:sparse_min_prob} with $K = AB$.
Indeed, solutions of this type of problem admit only finitely
many non-zero coefficients and often coincide with the sparsest
solution possible
\cite{donoho2004stablesparse,fuchs2005exactsparse,gribonval2007sparsenessmeasure}.

Unfortunately, the numerical solution of the above (non-smooth)
minimization problem is not straightforward. There is a vast amount of
literature dealing with efficient computational algorithms for
equivalent formulations of the problem
\cite{daubechies2007projgrad,figueiredo2007gradproj,
  turlach2005algorithmsleastsquaresL1,kim2007l1ls,osborne2000variableselection,efron2004lars,griesse2007ssnsparsity,elad2007subspaceoptimization}, both
in the infinite-dimensional setting as well as for finitely many
dimensions, but mostly for the finite-dimensional case.
An often-used, simple but apparently slow algorithm is the
iterative soft-thresholding (or thresholded Landweber) procedure 
which is known to
converge in the strong sense in infinite dimensions~\cite{daubechies2003iteratethresh}.
The algorithm is simple: it just needs an initial value $u^0$
and an operator with $\norm{K} < 1$. The iteration reads as
\[
    u^{n+1} = \mathbf{S}_{\alpha}\bigl( u^n - K^*(Ku^n - f) \bigr)
    \quad , \quad \bigl(\mathbf{S}_{\alpha}(w)\bigr)_k = \sign{w_k}
    \bigl[ \abs{w_k} - \alpha_k \bigr]_+.
\]
In practice it is important to know moreover 
convergence rates for the algorithms or at least an estimate for the
distance to a minimizer to evaluate the fidelity of the
outcome of the computations.
The convergence proofs in the infinite-dimensional case presented in
\cite{daubechies2003iteratethresh}, and for generalizations in
\cite{combettes2005signalrecovery}, however, 
do not imply a-priori estimates and do not inherently give
any rate of convergence, although, in many cases, linear convergence 
can be deduced quite easily from the fact that iterative thresholding 
converges strongly and from the special structure of the algorithm. 
To the best knowledge of the authors, \cite{bredies2008harditer} 
contains the first results about the convergence of 
iterative algorithms for linear inverse problems 
with sparsity constraints in infinite dimensions for which the 
convergence rate is inherent in the respective proof. There, an
iterative hard-thresholding procedure has been proposed for which, if
$K$ is injective, a convergence rate of $\mathcal{O}(n^{-1/2})$ could be
established. 

The main purpose of this paper is to develop a general and unified framework for the convergence analysis of algorithms for the problem~(\ref{eq:sparse_min_prob}) and related problems, especially for the iterative soft-thresholding algorithm.
We show that the iterative soft-thresholding algorithm converges linearly in almost every case and point out how to obtain a-priori estimates.
To this end, we formulate the iterative soft-thresholding as a generalized gradient projection method which leads to a new proof for the strong convergence which is independent of the proof given in~\cite{daubechies2003iteratethresh}.
The techniques used for our approach 
may shed new light on the known 
properties of the iterative soft-thresholding related methods. 

We distinguish two key properties which lead to linear convergence.
The first is called \emph{finite basis injectivity} (FBI) and is a
property of the operator $K$ only, while the second is called 
a \emph{strict sparsity pattern} of a solution of the minimization
problem~(\ref{eq:sparse_min_prob}).
\begin{Definition}
  \label{def:fbi}
  An operator $K:\lpspace{2}\to\HH_2$ mapping into a Hilbert
  space has the
  \emph{finite basis injectivity} property, if for all finite subsets
  $I\subset\NN$ the operator $K|_I$ is injective, i.e.~for
  all $u,v\in \lpspace{2}$ with $Ku = Kv$ and
  $u_k = v_k = 0$ for all $k \notin I$ it follows $u = v$.
\end{Definition}
\begin{Definition}
  \label{def:strict_sparsity_pattern}
  A solution $u^*$ of~(\ref{eq:sparse_min_prob}) possesses a
  \emph{strict sparsity pattern} if whenever $u_k^* = 0$ for some $k$
  there follows $\abs{K^*(Ku^* - f)}_k < \alpha_k$.
\end{Definition}

The main result can be summarized by the following:
\begin{Theorem}
  \label{thm:iter_thres_lin_conv}
  Let $K: \lpspace{2} \rightarrow \HH_2$, $K \neq 0$
  be a linear and continuous
  operator as well as $f \in \HH_2$.
  Consider the sequence $\seq{u^n}$ given by the iterative
  soft-thresholding procedure 
  \begin{equation}
    \label{eq:iter_soft_thres}
    u^{n+1} = \mathbf{S}_{s_n\alpha}\bigl( u^n - s_nK^*(Ku^n - f) \bigr)
    , \quad \bigl(\mathbf{S}_{s_n\alpha}(w)\bigr)_k = \sign{w_k}
    \bigl[ \abs{w_k} - s_n\alpha_k \bigr]_+
  \end{equation}
  with step size 
  \begin{equation}
    \label{eq:iter_soft_thres_step_size}
    0 < \underline{s} \leq s_n \leq \overline{s} < 2/\norm{K}^2 
  \end{equation}
  and a $u^0 \in \lpspace{2}$ such that $\sum_{k=1}^\infty \alpha_k
  \abs{u_k^0} < \infty$.
  Then, there is a minimizer $u^*$ such that $u^n \rightarrow u^*$ in
  $\lpspace{2}$.
  
  Moreover, suppose that either
  \begin{enumerate}
  \item $K$ possesses the FBI property, or
  \item $u^*$ possesses a strict sparsity pattern.
  \end{enumerate}
  Then, $u^n \rightarrow u^*$ with a linear rate, i.e.~there exists a
  $C > 0$ and a $0 \leq \lambda < 1$ such that $\norm{u^n - u^*} \leq
  C \lambda^n$.
\end{Theorem}

\begin{Remark}[(Examples for operators with the FBI property)]
  In the context of inverse problems with sparsity constraints, the FBI property is natural, since the operators $A$ are often injective.
  Prominent examples are the Radon transform \cite{maass1992interiorradon}, solution operators for partial differential equations, e.g.~in heat conduction problems~\cite{dahlke2003waveletgalerkin} or inverse boundary value problems like electrical impedance tomography~\cite{nachman1996invboundaryvalueproblem}.
  The combination with a synthesis operator $B$ for an orthonormal basis does not influence the injectivity.

  Moreover, the restriction to orthonormal bases can be relaxed.
  The results presented in this paper also hold if the system $\{\psi_k\}$ is a frame or even a dictionary---as long as the FBI property is fulfilled.
  This is for example the case for a frame which consists of two orthonormal bases where no element of one basis can be written as a finite linear combination of elements of the other.
  This is typically the case, e.g.~for a trigonometric basis and the Haar wavelet basis on a compact interval.
  One could speak of FBI frames or FBI dictionaries.
\end{Remark}

\begin{Remark}[(Strict sparsity pattern)]
  This condition can be interpreted as follows. We know that the
  weighted $\lpspace{1}$-regularization imposes sparsity on
  a solution $u^*$ in the sense that $u_k^* = 0$ for all but finitely
  many $k$, hence the name sparsity constraint. 
  For the remaining indices, the equations $(K^*Ku^*)_k = K^*f
  - \alpha_k \sign{u^*_k}$ are satisfied which corresponds to an
  approximate solution of the generally ill-posed equation $Ku = f$ in a
  certain way. Now the condition that the solutions of
  (\ref{eq:sparse_min_prob}) possess a strict sparsity pattern says
  that $u_k^* = 0$ for some index $k$ can occur only because of the
  sparsity constraint but never for the solution of the linear
  equation.
  We emphasize that Theorem~\ref{thm:iter_thres_lin_conv} states that
  whenever $\seq{u^n}$ converges to a solution $u^*$ with strict
  sparsity pattern, then the speed of convergence has to be linear
  for all bounded linear operators $K$. 
\end{Remark}

The proof of Theorem~\ref{thm:iter_thres_lin_conv} will be divided
into three sections. First, in Section \ref{sec:iter_thres_grad_proj},
we introduce a
framework in which iterative soft-thresholding according to
\eqref{eq:iter_soft_thres} can be interpreted as a generalized
gradient projection method. We derive descent properties for
generalized gradient methods and show under which conditions we can
obtain linear convergence in Section \ref{sec:conv_gen_grad_proj}.
We show in Section \ref{sec:conv_rates} 
that a Bregman-distance estimate for problems of the type
\eqref{eq:sparse_min_prob} gives a new convergence proof for the
iterative soft-thresholding.
In Section~\ref{sec:conv-relat-meth} we illustrate the 
broad range of applicability of the results with two more examples.
Finally, some conclusions about the implications of the
results are drawn in Section \ref{sec:conclusions}.

\section{Iterative soft-thresholding and a
  generalized gradient projection method}
\label{sec:iter_thres_grad_proj}

A common approach to solve smooth unconstrained minimization problems
are methods based on moving in the direction of steepest descent, i.e.~the 
negative gradient. 
In constrained optimization, the gradient is often projected back to the 
feasible set, yielding the well-known gradient projection algorithm
method \cite{goldstein1964convexprog,levitin1966constrainedmin,
  dunn1981gradientprojection}. 
In the following, a step of generalization is introduced: 
The method is extended to deal with
sums of smooth and non-smooth functionals, and covers in particular
constrained smooth minimization problems. The gain is that the iteration
\eqref{eq:iter_soft_thres} fits into this generalized framework.

Similar to the generalization performed in \cite{bredies2005gencondgrad}, 
its main idea is to replace the constraint by a general proper, convex and
lower semi-continuous functional $\Phi$ which leads, 
for the gradient projection method, to the successive 
application of the associated proximity operators, i.e.
\begin{equation}
  \label{eq:prox_operator}
  J_s: w\mapsto \argmin_{v\in\HH} \ \frac{\norm{v - w}^2}{2} + s\Phi(v)\ .
\end{equation}
The generalized gradient projection method for minimization problems 
of type
\begin{equation}
  \label{eq:min_prob_sum}
  \min_{u \in \HH} \ F(u) + \Phi(u)
\end{equation}
then read as follows.
\begin{Algorithm}
  \label{alg:gen_grad_proj}
  \mbox{}

  \begin{enumerate}
  \item Choose a $u^0 \in \HH$ with $\Phi(u^0) < \infty$ and set $n =
    0$.
  \item Compute the next iterate $u^{n+1}$ according to
    \[ u^{n+1} = J_{s_n} \bigl( u^n - s_n F'(u^n) \bigr) \ . \]
    where $s_n$ satisfies an appropriate step-size rule and 
    $J_s$ from~(\ref{eq:prox_operator}).
  \item Set $n := n + 1$ and continue with Step 2.
  \end{enumerate}
\end{Algorithm}
Note that the solutions of the minimization problem are exactly the fixed 
points of the algorithm.
Moreover, the case $\Phi= I_\Omega$, where $\Omega$ is a closed and convex 
constraint, yields the classical gradient projection method 
which is known to converge
provided that certain assumptions are
fulfilled and a suitable step-size rule has been chosen
\cite{dunn1981gradientprojection,demyanov1970approxoptim}.

In the following, we assume that $F$ is differentiable, $F'$ is Lipschitz 
continuous with constant $L$ and usually choose the step-sizes such that
\begin{equation}
  \label{eq:step_size_constraint}
  0 < \underline{s} \leq s_n \leq \overline{s} < 2/L. 
\end{equation}
Note that form the trivial case $L=0$ we agree that $2/L = \infty$.

\begin{Remark}[(Forward-backward splitting)] 
  \label{rem:forward_backward_splitting}
  The generalization of the gradient projection method leads to a
  special case of the so-called \emph{proximal forward-backward splitting
    method}
  which amounts to the iteration
  \[
  u^{n+1} = u^{n} + t_n \Bigl( J_{s_n}\bigl( u - s_n ( F'(u^n)
  + b^n) \bigr) + a^n - u^{n} \Bigr)
  \]
  where $t_n \in [0,1]$ and $\seq{a^n}, \seq{b^n}$ are
  absolutely summable sequences in $\HH$.
  In \cite{combettes2005signalrecovery}, it is shown that this method converges strongly
  to a minimizer under appropriate conditions. There exist, however,
  no general statements about convergence rates so far.
  Here, we restrict ourselves to the special case of the
  generalized gradient projection method. 

\end{Remark}

Finally, it is easy to see that 
the iterative soft-thresholding algorithm \eqref{eq:iter_soft_thres}
is a special case of this generalized gradient projection method
in case the functionals $F:\lpspace{2}\to\RR$ and $\Phi:\lpspace{2}\to
{]{-\infty},\infty]}$ are
chosen according to
\begin{equation}
  \label{eq:sparsity_functional_split}
  F(u) = \frac{\norm{Ku-f}^2}{2} \quad , \quad
  \Phi(u) = \begin{cases}
                \sum_{k=1}^\infty \alpha_k \abs{u_k}\ , & \text{ if the sum converges}\\
                \infty\ , & \text{else}
            \end{cases}
\end{equation}
where $K: \lpspace{2} \rightarrow \HH_2$ is
linear
and continuous between the Hilbert spaces $\lpspace{2}$ and $\HH_2$,
$f \in \HH_2$ and $\seq{\alpha_k}$ is sequence satisfying $\alpha_k
\geq \underline{\alpha} > 0$ for all $k$.

Here, $F'(u) = K^*(Ku - f)$, so in each iteration step of
Algorithm~\ref{alg:gen_grad_proj} we have to solve
\[
\min_{v \in \lpspace{2}} \ \frac{\norm{u^n - s_n K^*(Ku^n - f) -
    v}^2}{2} + s_n \sum_{k=1}^\infty \alpha_k \abs{v_k} 
\]
for which the solution is given by soft-thresholding, i.e.
\[
v = \mathbf{S}_{s_n\alpha}\bigl(u^n - s_n K^*(Ku^n - f) \bigr) \ ,
\]
with $\mathbf{S}_{s_n\alpha}$ according to~\eqref{eq:iter_soft_thres},
see \cite{daubechies2003iteratethresh}, for example.

Since the Lipschitz constant associated with $F'$ does not exceed
$\norm{K}^2$, this result can be summarized as follows:
\begin{Proposition}
  \label{prop:threshold_gen_grad_proj_equiv}
  Let $K:\lpspace{2}\to\HH_2$ be a bounded linear operator,
  $f\in\HH_2$ and $0 < \underline{\alpha} \leq \alpha_k$. 
  Let $F$ and $\Phi$ be chosen according to
  \eqref{eq:sparsity_functional_split}. Then Algorithm~\ref{alg:gen_grad_proj}
  with step-size $\seq{s_n}$ according to 
  \eqref{eq:iter_soft_thres_step_size}
  coincides with the iterative soft-thresholding procedure
  \eqref{eq:iter_soft_thres}.
\end{Proposition}
Here and in the following,
we also agree to set $2/\norm{K}^2 = \infty$ in
\eqref{eq:iter_soft_thres_step_size} for the trivial case $K=0$.

\section{Convergence of the generalized gradient 
  projection method}
\label{sec:conv_gen_grad_proj}

In the following, conditions which ensure convergence of the
generalized gradient projection method are derived.
The key is the descent of
the functional $F + \Phi$ in each iteration step. 
The following lemma states some basic properties of one
iteration.
\begin{Lemma}
  \label{lem:basic_descent_properties}
  Let $F$ be differentiable with $F'$ Lipschitz continuous with
  associated constant $L$ and $\Phi$
  be proper, convex and lower semi-continuous.
  Set $v = J_s\bigl( u - s F'(u) \bigr)$ as in (\ref{eq:prox_operator}) for some $s > 0$ 
  and denote by
  \begin{equation}
    \label{eq:bregman_like_distance}
    D_s(u) = \Phi(u) - \Phi(v) + \scp{F'(u)}{u - v}
  \end{equation}
  Then it holds:
    \begin{gather}
      \label{eq:dist_scp}
       \forall w\in\HH:\ \Phi(w) - \Phi(v) + \scp{F'(u)}{w - v}
       \geq \frac{\scp{u - v}{w - v}}{s}.\\
      \label{eq:dist_norm2}
      D_s(u) \geq \frac{\norm{v-u}^2}{s}\\
    \label{eq:lip_f_decr}
      (F+\Phi)(v) \leq (F+\Phi)(u) - \Bigl(1 - \frac{sL}{2}\Bigr)D_s(u).
    \end{gather}
\end{Lemma}

\begin{Proof}
  Since $v$ solves the problem
  \[
  \min_{v \in \HH} \ \frac{\norm{v - u + sF'(u)}^2}{2} + s \Phi(v)
  \]
  it immediately follows that the subdifferential inclusion $u - sF'(u) - v
  \in s\subgrad \Phi(v)$ is satisfied, see
  \cite{ekelandtemam1976convex,rockafellar1998variation} for an
  introduction to
  convex analysis and subdifferential calculus. This can be rewritten
  to
  \[
  \scp{u - sF'(u) - v}{w - v} \leq s \bigl( \Phi(w) 
  - \Phi(v) \bigr) \quad \text{for all} \ w \in \HH \ ,
  \]
  while rearranging and dividing by $s$ proves the
  inequality~(\ref{eq:dist_scp}).
  The inequality~(\ref{eq:dist_norm2}) follows by 
  setting $w = u$ in (\ref{eq:dist_scp}).
  
  To show inequality~(\ref{eq:lip_f_decr}), we observe
  \begin{multline*}
    (F+\Phi)(v) - (F+\Phi)(u) +D_s(u) = F(v) - F(u) + \scp{F'(u)}{u-v}\\
     = \int_0^1 \scp{F'\bigl( u + t(v-u) \bigr) - F'(u)}{v - u}
     \dd{t} \ .
  \end{multline*}
  Using the Cauchy-Schwarz inequality and the Lipschitz continuity we obtain
  \begin{align*}
    (F+\Phi)(v) - (F+\Phi)(u) +D_s(u) & \leq \int_0^1 t
    L\norm{v-u}^2\dd{t}
    = \tfrac{L}{2}\norm{v-u}^2.
  \end{align*}
    Finally, applying the estimate~(\ref{eq:dist_norm2})
    leads to~(\ref{eq:lip_f_decr}).
\end{Proof}

\begin{Remark}[(A weaker step-size condition)]
  \label{rem:step_size_descent}
  If the step-size in the generalized gradient projection method 
  is chosen such that
  $s_n \leq \overline{s} < 2/L$, then
  we can conclude from~\eqref{eq:lip_f_decr} that
  \begin{equation}
    \label{eq:rest_distance_est}
  (F+\Phi)(u^{n+1}) \leq (F+\Phi)(u^n) - \delta D_{s_n}(u^n)
  \end{equation}
  where $\delta = 1 - \frac{\overline{s}L}{2}$.
  Of course, the constraint on the step size is only sufficient to
  guarantee such a decrease.
  A weaker condition is the following:
  \begin{equation}
    \label{eq:weak_step_size_cond}
    \Bigabs{\int_0^1 \scp{F'\bigl( u^n + t(u^{n+1}-u^n) \bigr) -
        F'(u^n)}{u^{n+1} - u^n} 
      \dd{t} }  \leq (1 - \delta) D_{s_n}(u^n)
  \end{equation}
  for some $\delta > 0$. Regarding the proof of Lemma
  \ref{lem:basic_descent_properties}, it is easy to see that this
  condition also leads to the estimate
  \eqref{eq:rest_distance_est}. Unfortunately,
  \eqref{eq:weak_step_size_cond} can only
  be verified a-posteriori, i.e.~with the knowledge of the next
  iterate $u^{n+1}$. So one has to guess an $s_n$ and check if
  \eqref{eq:weak_step_size_cond} is satisfied, otherwise a different
  $s_n$ has to be chosen. In practice, this means that one iteration
  step is lost and consequently more computation time is needed,
  reducing the advantages of a more flexible step size.
\end{Remark}

While the descent property \eqref{eq:rest_distance_est} can be proven
without convexity assumptions on $F$, we need such a property to
estimate the distance of the functional values 
to the global minimum of $F+\Phi$ in the
following. 
We introduce for any sequence $\seq{u^n} \subset \HH$ according to
Algorithm \ref{alg:gen_grad_proj} the values
\begin{equation}
  \label{eq:r_n}
  r_n = (F+ \Phi)(u^n) - \Bigl( \min_{u \in \HH} \ (F+ \Phi)(u)  \Bigr) \ .
\end{equation}

\begin{Proposition}
  \label{prop:descent_rate}
  Let $F$ be convex and continuously differentiable with Lip\-schitz
  continuous derivative.
  Let $\seq{u^n}$ be a sequence generated by Algorithm~\ref{alg:gen_grad_proj}
  such that 
  the step-sizes are bounded from below, i.e.~$s_n \geq \underline{s} > 0$,
  and that we have
  \[
  (F+\Phi)(u^{n+1}) \leq (F+\Phi)(u^n) - \delta D_{s_n}(u^n)
  \]
  for a $\delta > 0$ with $D_{s_n}(u^n)$ according to
  \eqref{eq:bregman_like_distance}.

  \begin{enumerate}
  \item  If $F + \Phi$ is coercive, then
    the values $r_n$ according to \eqref{eq:r_n} satisfy
    $r_n \rightarrow 0$ with rate $\kronO(n^{-1})$,
    i.e.~there exists a $C > 0$ such that
    \[
    r_n \leq C n^{-1} \ .
    \]
  \item If for a minimizer $u^*$ and some $c > 0$ 
    the values $r_n$ from \eqref{eq:r_n} satisfy 
    \begin{equation}
      \label{eq:rest_norm_est}
      \norm{u^n - u^*}^2
      \leq c r_n \ ,
    \end{equation}
    
    then $\seq{r_n}$ vanishes exponentially and
    $\seq{u^n}$ converges linearly to $u^*$, i.e.~there exists a $C > 0$
    and a $\lambda \in [0,1[$ such that
    \[
    \norm{u^n - u^*} \leq C \lambda^n \ .
    \]
  \end{enumerate}
\end{Proposition}

\begin{Proof}
  We first prove an estimate for $r_n$ and then treat the cases
  separately. For this purpose,
  pick an optimal $u^* \in \HH$ and 
  observe that the decrease in each iteration step can be estimated
  by
  \[
  r_n - r_{n+1} = (F + \Phi)(u^n) - (F + \Phi)(u^{n+1}) \geq \delta
  D_{s_n}(u^n) \ ,
  \]
  according to the assumptions.
  Note that 
  $D_{s_n}(u^n) \geq 0$ by \eqref{eq:dist_norm2}, so $\seq{r_n}$ is
  non-increasing. 
  
  Use the convexity of $F$ to deduce
  \begin{align*}
    r_n 
    &\leq \Phi(u^n) - \Phi(u^*) + \scp{F'(u^n)}{u^n - u^*} \displaybreak[2]\\
    &= D_{s_n}(u^n) + \scp{F'(u^n)}{u^{n+1} - u^*} + \Phi(u^{n+1}) -
    \Phi(u^*) \\
    &\leq D_{s_n}(u^n) + \frac{\scp{u^n - u^{n+1}}{u^{n+1} - u^*}}{s_n} \\
    & \leq D_{s_n}(u^n) + \frac{\norm{u^{n+1} - u^*}}{\sqrt{s_n}}
    \sqrt{D_{s_n}(u^n)}
  \end{align*}
  by applying the Cauchy-Schwarz inequality 
  as well as \eqref{eq:dist_scp} and
  \eqref{eq:dist_norm2}.
  With the above estimate on $r_n - r_{n+1}$ and $0 < \underline{s} <
  s_n$ we get
  \begin{equation}
    \label{eq:rest_proof_est}
    \delta r_n \leq (r_n - r_{n+1}) + \frac{\sqrt{\delta}\norm{u^{n+1} -
        u^*}}{\sqrt{\underline{s}}} \sqrt{r_n - r_{n+1}} \ .
  \end{equation}
  
  We now turn to prove the first statement of the proposition.
  Assume that $F + \Phi$ is coercive, so from the fact that 
  $\seq{r_n}$ is non-increasing follows that
  $\norm{u^n - u^*}$ has to be bounded by a $C_1 > 0$. Furthermore,
  $0 \leq r_n - r_{n+1} \leq r_0 < \infty$, implying
  \[
  \delta r_n \leq \bigl( \sqrt{r_0} + \sqrt{\delta \underline{s}^{-1} }
  C_1 \bigr) \sqrt{r_n - r_{n+1}}
  \]
  and consequently
  \[
  q r_n^2 \leq r_n - r_{n+1} \quad , \quad q =
  \Bigl( \frac{\delta}{\sqrt{r_0} + \sqrt{\delta \underline{s}^{-1} }
    C_1} \Bigr)^2 > 0 \ .
  \]
  Standard arguments then give the rate $r_n = \kronO(n^{-1})$, we
  repeat them here for convenience. The above estimate on $r_n -
  r_{n+1}$ as well the property that $\seq{r_n}$ is non-increasing
  yields
  \[
  \frac{1}{r_{n+1}} - \frac{1}{r_n} = \frac{r_n - r_{n+1}}{r_nr_{n+1}}
  \geq q \frac{r_n^2}{r_n r_{n+1}} \geq q
  \]
  which, summed up, leads to 
  \[
  \frac{1}{r_n} - \frac{1}{r_0} = \sum_{i=0}^{n-1} \frac{1}{r_{i+1}} -
  \frac{1}{r_i} \geq nq \quad \Rightarrow \quad r_n^{-1} \geq nq +
  r_0^{-1}
  \]
  and consequently, since $q >0$, to the desired rate 
  $r_n \leq (nq + r_0^{-1})^{-1} \leq C n^{-1}$.  
  
  Regarding the second statement, assume that there is a $c > 0$ such
  that $\norm{u^n - u^*}^2 \leq c r_n$ for some optimal $u^*$ and each
  $n$. Starting again at \eqref{eq:rest_proof_est} and applying
  Young's inequality yields, for each $\varepsilon > 0$,
  \[
  \delta r_n \leq (r_n - r_{n+1}) + \frac{\delta\varepsilon \norm{u^{n+1} -
      u^*}^2}{2\underline{s}} + \frac{r_n - r_{n+1}}{2\varepsilon}
  \ .
  \]
  Choosing $\varepsilon = \underline{s} c^{-1}$ and exploiting the 
  assumption $\norm{u^{n+1} - u^*}^2 \leq cr_{n+1}$ 
  as well as the fact $r_{n+1} \leq r_n$ 
  then imply 
  \[
  \delta r_n \leq (r_n - r_{n+1}) + \frac{\delta}{2} r_n + 
  \frac{r_n - r_{n+1}}{2 \underline{s} c^{-1}} 
  \quad \Rightarrow \quad
  r_n - r_{n+1} \geq \frac{\delta \underline{s} c^{-1}}{2
    \underline{s} c^{-1} + 1} r_n
  \]
  which in turn establishes the exponential decay rate
  \begin{equation}
    \label{eq:rest_exp_decay}
    r_{n+1} \leq  \Bigl( 1 - \frac{\delta \underline{s}
      c^{-1}}{2\underline{s} c^{-1} + 1} \Bigr) r_n \leq
    \lambda^2 r_n
    \ ,\ \text{ with }
    \lambda =  \Bigl( 1 - \frac{\delta \underline{s}
      c^{-1}}{2\underline{s} c^{-1} + 1} \Bigr)^{1/2} \in [0,1[.
  \end{equation}
  Using $\norm{u^n - u^*}^2 \leq c r_n$ again finishes
  the proof:
  \[
  \norm{u^n - u^*} \leq (cr_n)^{1/2} \leq (cr_0)^{1/2} \lambda^n \ .  
  \]
\end{Proof}

Proposition \ref{prop:descent_rate} tells us that we only have
to establish \eqref{eq:rest_norm_est} to obtain strong convergence
with linear convergence rate. This can be done with determining
how fast the functionals $F$ and $\Phi$ vanish at some minimizer. 
This can be made precise by introducing the following notions which
also turn out to be the essential ingredients to show
\eqref{eq:rest_norm_est}:
First, define for a minimizer $u^* \in \HH$ the functional
\begin{equation}
  \label{eq:strong_conv_est}
  R(v) = \scp{F'(u^*)}{v - u^*} + \Phi(v) - \Phi(u^*) \ .
\end{equation}
Note that if the subgradient of $\Phi$ in $u^*$ is unique, $R$ is the
Bregman distance of $\Phi$ in $u^*$, a notion which is extensively
used in the analysis of descent algorithms \cite{schoepfer2006illposedbanach,bauschke2003bregman}.
Moreover, we make use of the remainder of the Taylor expansion of $F$,
\begin{equation}
  \label{eq:taylor_rest}
  T(v) = F(v) - F(u^*) - \scp{F'(u^*)}{v - u^*} \ .
\end{equation}

\begin{figure}[tb]
  \centering
  \begin{pgfpicture}{0cm}{0cm}{8cm}{6cm}
    \pgfputat{\pgfxy(0,0)}{\pgfbox[left,bottom]{\includegraphics[width=8cm]{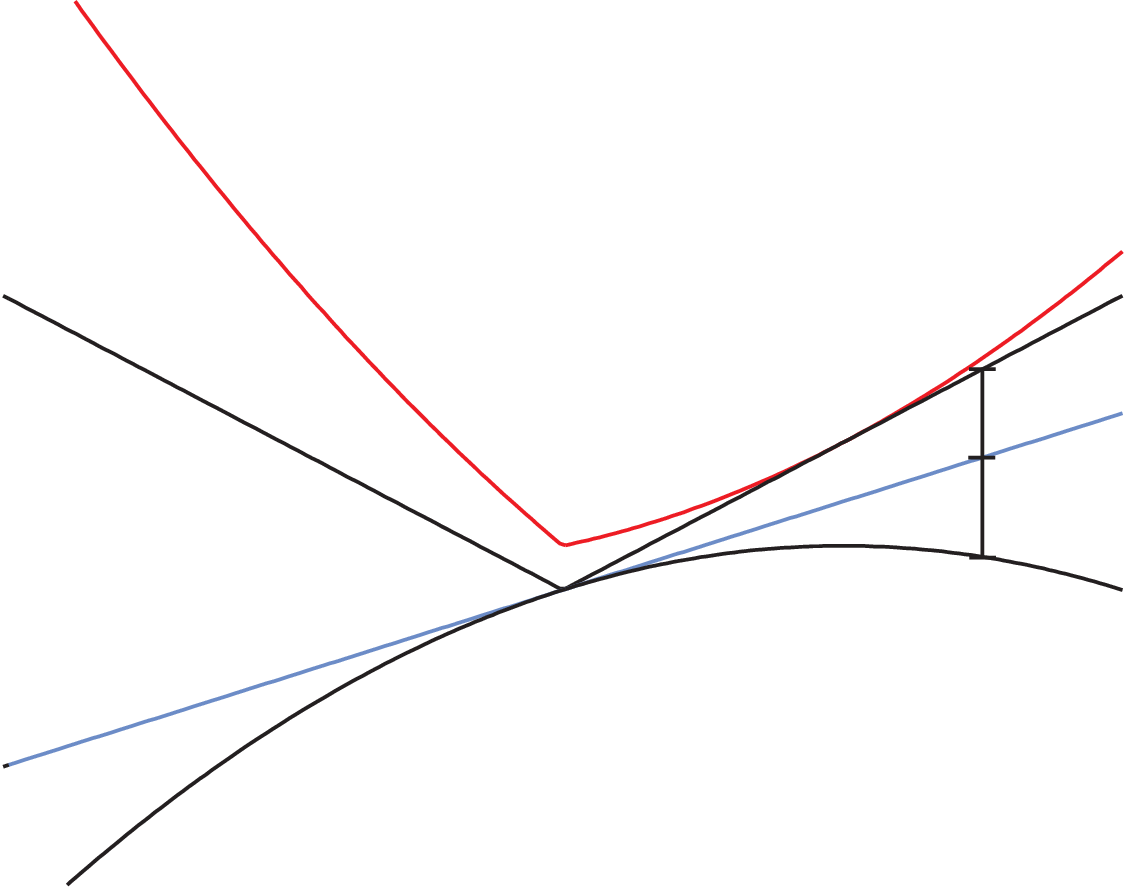}}}
    \pgfputat{\pgfxy(3.6,4)}{\pgfbox[center,center]{$(F+\Phi)(v)$}}
    \pgfputat{\pgfxy(1.7,3.7)}{\pgfbox[center,center]{$\Phi(v)$}}
    \pgfputat{\pgfxy(3.1,0.5)}{\pgfbox[center,center]{$-F(v)+
        (F+\Phi)(u^*)$}}
    \pgfputat{\pgfxy(0.8,2)}{\pgfbox[center,center]{$\Phi(u^*)-\scp{F'(u^*)}{v
          - u^*}$}}
    \pgfputat{\pgfxy(7.2,3.45)}{\pgfbox[center,center]{$R$}}
    \pgfputat{\pgfxy(7.2,2.7)}{\pgfbox[center,center]{$T$}}
  \end{pgfpicture}
  
  \caption{
    \label{fig:bregman-taylor-distance}
    Illustration of the Bregman-like distance $R$ and the
    Taylor distance $T$ for a convex $\Phi$ and a smooth $F$.
    Note that for the optimal value $u^*$ it holds $-F'(u^*)\in
    \subgrad\Phi(u^*)$.}
\end{figure}

\begin{Remark}[(On the Bregman distance)]
  \label{rem:bregman_sufficient}
  In many cases the Bregman-like distance
  $R$ is enough to estimate the descent properties,
  see \cite{bredies2008harditer,schoepfer2006illposedbanach}.
  For example, in
  case that $\Phi$ is the $p$-th power of a 
  norm of a $2$-convex Banach space $X$,
  i.e.~$\Phi(u) = \norm[X]{u}^p$ with $p \in {]{1,2}]}$,
  which is moreover continuously embedded in $\HH$, one can show that
  \begin{equation*}
    \norm[X]{v-u^*}^2 \leq C_1 R(v)
  \end{equation*}
  holds on each bounded set of $X$, see \cite{xu1991charinequalities}.
  Consequently, with $j_p = \subgrad \frac1p
  \norm[X]{\placeholder}^p$ denoting the duality mapping with gauge $t
  \mapsto t^{p-1}$,
  \begin{align*}
    \norm{v-u^*}^2 &\leq C_2 \norm[X]{v - u^*}^2 \\
    &\leq C_1 C_2 \bigl( \norm[X]{v}^p
    - \norm[X]{u^*}^p - p \scp{j_p(u^*)}{v - u^*} \bigr) = c R(v) 
  \end{align*}
  observing that $R$ is in this case the Bregman distance.
  Often, Tikhonov functionals for inverse problems admit such a
  structure, e.g.
  \[
  \min_{u \in \lpspace{2}} \ \frac{\norm{Ku - f}^2}{2\alpha} +
  \sum_{k=1}^\infty \abs{u_k}^p \ ,
  \]
  a regularization which is also topic in
  \cite{daubechies2003iteratethresh}.
  As one can see in complete analogy to Proposition
  \ref{prop:threshold_gen_grad_proj_equiv}, the generalized gradient
  projection method also amounts to the iteration proposed there, so
  as a by-product and after verifying that the prerequisites of
  Proposition \ref{prop:descent_rate} indeed hold, 
  one immediately gets a linearly convergent method.
  
  However, in the case that $\Phi$ is not sufficiently convex, the
  Bregman distance alone is not sufficient to obtain
  the required estimate on the $r_n$.
  This is the case for $F$ and $\Phi$ according
  \eqref{eq:sparsity_functional_split}.
  In this situation we also have to take the ``Taylor
  distance'' $T$ into account.
  Figure~\ref{fig:bregman-taylor-distance} shows an illustration of the values $R$ and $T$.
  One could say that the Bregman distance measures the error corresponding to the $\Phi$ part while the Taylor distance does the same for the $F$ part.
\end{Remark}

The functionals $R$ and $T$ possess the following properties:
\begin{Lemma}
  \label{lem:rest_estimate}
  Consider the problem \eqref{eq:min_prob_sum} where $F$ is
  convex, differentiable and $\Phi$ is proper, convex and lower
  semi-continuous. If $u^* \in \HH$ is a solution of
  \eqref{eq:min_prob_sum} and $v \in \HH$ is arbitrary, then the
  functionals $R$ and $T$ according to \eqref{eq:strong_conv_est} and
  \eqref{eq:taylor_rest}, respectively, are non-negative and satisfy
  \[
  R(v) + T(v) = (F+\Phi)(v)-(F+\Phi)(u^*) \ .
  \]
\end{Lemma}

\begin{Proof}
  The identity is obvious from the definition of $R$ and $T$.
  For the non-negativity of $R$, 
  note that since $u^*$ is a solution, it holds that $-F'(u^*) \in
  \subgrad \Phi(u^*)$. Hence, the subgradient inequality reads as
  \[
  \Phi(u^*) - \scp{F'(u^*)}{v - u^*} \leq \Phi(v) \quad \Rightarrow
  \quad R(v) \geq 0 \ .
  \]
  Likewise, the property $T(v) \geq 0$ is a consequence of the
  convexity of $F$.
\end{Proof}

Now it follows immediately that $R(v) = T(v) = 0$ whenever $v$ is a
minimizer. To conclude this section, 
the main statement about the convergence of the
generalized gradient projection method reads as:
\begin{Theorem}
  \label{thm:gen_grad_proj_lin_conv}
  Let $F$ be a convex, differentiable functional with
  Lipschitz-continuous derivative (with associated Lipschitz constant
  $L$), $\Phi$ be proper, convex and lower semi-continuous
  and $\seq{u^n}$ be a sequence generated by Algorithm
  \ref{alg:gen_grad_proj} with step-size according to
  \eqref{eq:step_size_constraint}. Moreover, suppose that $u^* \in
  \HH$ is a
  solution of the minimization problem \eqref{eq:min_prob_sum}.
  
  If, for each $M \in \RR$ there exists a constant $c(M) > 0$ such
  that
  \begin{equation}
    \label{eq:rest_local_estimate}
    \norm{v - u^*}^2 \leq c(M) \bigl(R(v) + T(v) \bigr)
  \end{equation}
  for each $v$ satisfying $(F+\Phi)(v) \leq M$ and $R(v)$ and $T(v)$ defined by (\ref{eq:strong_conv_est}) and (\ref{eq:taylor_rest}), respectively, then $\seq{u^n}$ converges
  linearly to the unique minimizer $u^*$.
\end{Theorem}

\begin{Proof}
  A step-size chosen according to \eqref{eq:step_size_constraint}
  implies, by Lemma \ref{lem:basic_descent_properties}, the descent
  property \eqref{eq:rest_distance_est} with $\delta = 1 -
  \overline{s} L/2$. 
  In particular, from \eqref{eq:rest_distance_est} follows that $\seq{r_n}$
  is non-increasing (also remember \eqref{eq:dist_norm2} means in
  particular that $D_{s_n}(u^n) \geq 0$).
  Now choose $M = (F+\Phi)(u^0) < \infty$ for which, by
  assumption, a $c > 0$ exists such that
  \begin{equation*}
    \norm{u^n - u^*}^2 \leq c\bigl(R(u^n) + T(u^n) \bigr) =
    c r_n.
  \end{equation*}
  Hence, the prerequisites for Proposition \ref{prop:descent_rate} are
  fulfilled and consequently, $u^n \rightarrow u^*$ with a linear
  rate. Finally, the minimizer has to be unique: If $u^{**}$ is also a
  minimizer, then $u^{**}$ plugged into
  \eqref{eq:rest_local_estimate} gives
  $\norm{u^{**} - u^*}^2 = 0$ and consequently $u^* = u^{**}$.
\end{Proof}

\section{Convergence rates for the iterative soft-thresholding
  method}
\label{sec:conv_rates}

We now turn to the proof of the main result, Theorem
\ref{thm:iter_thres_lin_conv}, which collects the results of Sections
\ref{sec:iter_thres_grad_proj} and \ref{sec:conv_gen_grad_proj}.
Within this section, we consider the regularized inverse 
problem~\eqref{eq:sparse_min_prob} under the prerequisites of 
Proposition~\ref{prop:threshold_gen_grad_proj_equiv}.
It is known that at least one minimizer for~\eqref{eq:sparse_min_prob}
exists~\cite{daubechies2003iteratethresh}.
%

We have already seen in Proposition
\ref{prop:threshold_gen_grad_proj_equiv} that 
the iterative thresholding procedure \eqref{eq:iter_soft_thres}
is equivalent to a generalized gradient projection method.
Our aim is, on the one hand, to apply
Proposition~\ref{prop:descent_rate} in order to get strong convergence
from the descent rate $\kronO(n^{-1})$.
On the other hand, we will show the applicability of
Theorem \ref{thm:gen_grad_proj_lin_conv} for $K$ possessing the FBI
property which implies the desired convergence speed.
Observe that $F$ and $\Phi$ meet the requirements of 
Theorem~\ref{thm:gen_grad_proj_lin_conv} and that the step-size 
rule~\eqref{eq:iter_soft_thres_step_size} 
immediately implies~\eqref{eq:step_size_constraint}, so we only have to 
verify~\eqref{eq:rest_local_estimate}.
This will be done, among a Bregman-distance estimate,
in the following lemma, which is also serving as the crucial
prerequisite for showing convergence.

\begin{Lemma}
  \label{lem:bregman_inf_dim_est}
  For each minimizer $u^*$ of~\eqref{eq:sparse_min_prob} and each $M
  \in \RR$, there exists a $c_1(M,u^*)$ and a subspace $U \subset
  \lpspace{2}$ with finite-dimensional complement such that for the
  Bregman-like distance \eqref{eq:strong_conv_est} it holds that
  \begin{equation}
    \label{eq:bregman_rest_sparsity_est}
    R(v) \geq c_1(M,u^*) \norm{P_U(v - u^*)}^2
  \end{equation}
  whenever $(F + \Phi)(v)
  \leq M$ with $F$ and $\Phi$ defined
  by~(\ref{eq:sparsity_functional_split}).

  If $K$ moreover satisfies the FBI property, there is a $c_2(M,u^*,K)
  > 0$ such that, whenever $(F+\Phi)(v) \leq M$, the associated
  Bregman-Taylor distance according
  to~\eqref{eq:strong_conv_est} and \eqref{eq:taylor_rest} satisfies
  \begin{equation}
    \label{eq:bregman_taylor_rest_sparsity_est}
    R(v) + T(v) \geq c_2(M,u^*,K) \norm{v - u^*}^2 \ .
   \end{equation}
\end{Lemma}

\begin{Proof}
  Let $u^*$ be a minimizer of~\eqref{eq:sparse_min_prob} and
  assume that $v \in \lpspace{2}$ satisfies $\Phi(v) \leq M$ for a $M
  \geq 0$. Then,
  \begin{equation}
    \label{eq:sparsity_bregman_dist}
    R(v) = \sum_{k=1}^\infty \alpha_k \abs{v_k} - \sum_{k=1}^\infty \alpha_k
    \abs{u_k^*} + \sum_{k=1}^\infty w_k^*(v_k - u_k^*)
  \end{equation}
  where $w^* = -F'(u^*) = -K^*(Ku^* - f)$. Now since $w^* \in \subgrad
  \Phi(u^*)$ we have $w_k^* \in \alpha_k \sgn(u_k^*)$ for each $k$
  (note that $\subgrad \abs{\placeholder} = \sgn(\placeholder)$ with
  $\sgn(0) = [-1,1]$), meaning that
  \[
  \alpha_k \bigl( \abs{v_k} - \abs{u_k^*} \bigr) + w_k^*(v_k - u_k^*)
  \geq 0
  \]
  for each $k$. Denote by $I = \set{k \geq 1}{\abs{w_k^*} = \alpha_k}$
  which has to be finite since $w^* \in \lpspace{2}$ implies
  \[
  \infty > \sum_{k \in I} \abs{w_k^*}^2 = \sum_{k \in I} \alpha_k^2 \geq
  \sum_{k \in I} \underline{\alpha}^2 = \abs{I} \underline{\alpha}^2 \ .
  \]
  Moreover, $w_k^* \rightarrow 0$ as $k \rightarrow \infty$, so there
  has to be a $\rho < 1$ such that $\abs{w_k^*}/\alpha_k \leq \rho$ for
  each $k \in \NN \without I$. Also, if $k \in \NN \without I$, then 
  $\abs{w_k^*} \leq \rho\alpha_k$ which means in particular that $u_k^*
  = 0$ since the opposite contradicts $w_k^* \in \alpha_k \sgn(u_k^*)$.
  So, one can estimate \eqref{eq:sparsity_bregman_dist}:
  \begin{align*}
    R(v) 
    &\geq \sum_{k \notin I} \alpha_k \abs{v_k} + w_k^*v_k 
    \geq \sum_{k \notin I} \alpha_k (1 - \rho) \abs{v_k} 
    \displaybreak[2]
    \\
    &\geq
    (1-\rho)\underline{\alpha} \sum_{k \notin I} \abs{v_k - u_k^*}
    \geq (1-\rho)\underline{\alpha} \Bigl( \sum_{k \notin I} \abs{v_k
      - u_k^*}^2 \Bigr)^{1/2}
  \end{align*}
  using the fact that one can estimate the $\lpspace{2}$-sequence norm
  with the $\lpspace{1}$-sequence norm, see \cite{bredies2008harditer}
  for example. 
  With $U = \set{v \in
    \lpspace{2}}{v_k = 0 \ \text{for} \ k \in I}$, the above also reads
  as $R(v) \geq (1-\rho)\underline{\alpha} \norm{P_U(v - u^*)}$ with
  $P_U$ being the orthogonal projection onto $U$ in $\lpspace{2}$. 
  
  Next, observe that $\underline{\alpha} \norm{P_Uv}\leq
  \underline{\alpha} \norm[1]{v} \leq \Phi(v) \leq M+1$, hence we have
  $\norm{P_Uv}^{-1} \geq
  \underline{\alpha}/(M + 1)$. Consequently,
  \[
  R(v) \geq \frac{(1-\rho) \underline{\alpha}^2}{M+1} \norm{P_U(v -
    u^*)}^2 = c_1(M,u^*) \norm{P_U(v - u^*)}^2 
  \]
  which corresponds to the estimate
  \eqref{eq:bregman_rest_sparsity_est}.
  Finally, $\Phi(v) \leq M$ whenever $(F+\Phi)(v) \leq M$ and there is
  no $v$ such that $(F+\Phi)(v) < 0$. Hence, for each $M \in \RR$
  there is a constant for which \eqref{eq:bregman_rest_sparsity_est}
  holds whenever $(F + \Phi)(v) \leq M$ which is the desired
  statement for $R$.

  To prove~\eqref{eq:bregman_taylor_rest_sparsity_est}, suppose $K$
  possesses the FBI property. Recall that $T(v)$ can be expressed by
  \begin{equation}
    T(v) = F(v) - F(u^*) - \scp{F'(u^*)}{v - u^*} = \frac{\norm{K(v - u^*)}^2}{2} \ .
    \label{eq:taylor_rest_sparsity_est}
  \end{equation}
  The claim now is that there is a $C(M,u^*,K)$ such that
  \[
  \norm{u}^2 \leq C(M, u^*,K) \bigl( c_1(M,u^*) \norm{P_Uu}^2
  + \tfrac12\norm{Ku}^2 \bigr)
  \]
  for each $u \in \lpspace{2}$. 
  We will derive this constant directly. First split $u = P_Uu +
  P_{U^\perp}u$, so we can estimate, with the help of the inequalities
  of Cauchy-Schwarz and Young ($ab \leq a^2/4 + b^2$ for $a,b \geq 0$),
  \begin{align*}
    \frac{\norm{Ku}^2}{2} &= \frac{\norm{KP_{U^\perp}u}^2}{2} +
    \scp{KP_{U^\perp}u}{KP_Uu} + \frac{\norm{KP_Uu}^2}{2}\\
    &\geq \frac{\norm{KP_{U^\perp}u}^2}{4} - \frac{\norm{KP_Uu}^2}{2} 
    \geq \frac{\norm{KP_{U^\perp}u}^2}{4} - \frac{\norm{K}^2}{2}
    \norm{P_Uu}^2 \ .
  \end{align*}
  Since $K$ fulfills the FBI property, the
  operator restricted to $U^\perp$ is injective on the
  finite-dimensional space $U^\perp$, so there exists a $\bar c(U,K) >
  0$ such that $\bar c(U,K) \norm{P_{U^\perp}u}^2 \leq \norm{KP_{U^\perp}
    u}^2$ for all $u \in \lpspace{2}$. Hence,
  \[
  \norm{P_{U^\perp}u}^2 \leq 4 \bar c(U,K)^{-1} \bigl( \tfrac12 \norm{K}^2
  \norm{P_Uu}^2 + \tfrac12 \norm{Ku}^2 \bigr)
  \]
  and consequently
  \begin{align*}
    \norm{u}^2 &= \norm{P_{U^\perp}u}^2 + \norm{P_Uu}^2 \\
    &\leq 4\bar c(U,K)^{-1} \bigl( \bigl(\tfrac12 \norm{K}^2 + \tfrac14
    \bar c(U,K) \bigr) \norm{P_Uu}^2 + \tfrac12 \norm{Ku}^2 \bigr) \\
    & \leq \frac{2 \norm{K}^2 + \bar c(U,K) + 4c_1(M,u^*)}{\bar c(U,K)
      c_1(M,u^*)} 
    \bigl( c_1(M,u^*) \norm{P_Uu}^2 + \tfrac12 \norm{Ku}^2 \bigr) 
  \end{align*}
  giving a constant $c(M,u^*,K) > 0$ since $U$ depends on $u^*$.
  
  This finally yields the statement
  \begin{align*}
    \norm{v - u^*}^2 &\leq c(M, u^* ,K) \bigl( c_1(M, u^*)
    \norm{P_U(v - u^*)}^2 + \tfrac12 \norm{K(v - u^*)}^2  \bigr) \\
    &\leq 
    c(M,u^*,K) \bigl( R(v) + T(v) \bigr) \ ,
  \end{align*}
  consequently, \eqref{eq:bregman_taylor_rest_sparsity_est} holds for
  $c_2(M,u^*,K) = c(M,u^*,K)^{-1}$.
\end{Proof}

In the following, we will see that 
the estimate~\eqref{eq:bregman_rest_sparsity_est} considered in
$R(u^n)$ already leads to strong convergence of the iterative
soft-thresholding procedure. Nevertheless, we utilize
\eqref{eq:bregman_taylor_rest_sparsity_est} later, when the linear
convergence result will be proven.

\begin{Lemma}
  \label{lem:iter_thres_strong_conv}
  Let $K: \lpspace{2} \rightarrow \HH_2$ be a linear and continuous
  operator as well as $f \in \HH_2$. Consider the sequence $\seq{u^n}$
  which is generated by the iterative soft-thresholding procedure
  \eqref{eq:iter_soft_thres} with step-sizes $\seq{s_n}$
  according to \eqref{eq:iter_soft_thres_step_size}.
  Then, $\seq{u^n}$ converges to a minimizer in the strong sense.
\end{Lemma}

\begin{Proof}  
  Since the Lipschitz constant for $F'$ satisfies $L \leq \norm{K}^2$,
  the step sizes are fulfilling
  \eqref{eq:step_size_constraint} which
  implies, by Lemma \ref{lem:basic_descent_properties}, the descent
  property \eqref{eq:rest_distance_est} with $\delta = 1 -
  \overline{s}\norm{K}^2/2$.
  This means in particular that the associated functional distances
  $\seq{r_n}$ are non-increasing
  (since \eqref{eq:dist_norm2} in particular gives that $D_{s_n}(u^n)
  \geq 0$). Moreover, the descent result in Proposition
  \ref{prop:descent_rate} yields that the iterates $\seq{u^n}$ satisfy
  $R(u^n) \leq r_n \leq \kronO(n^{-1})$. Since $(F+\Phi)(u^n) 
  \leq (F+\Phi)(u^0) = M$, 
  we can apply Lemma \ref{lem:bregman_inf_dim_est} and
  the estimate \eqref{eq:bregman_rest_sparsity_est}
  leads to strong convergence of $\seq{P_Uu^n}$, i.e.~$P_U u^n \rightarrow P_U u^*$.
  
  Next, consider the complement parts $\seq{P_{U^\perp}u^n}$ in the
  finite-dimensional space $U^{\perp}$. Since $\norm{P_{U^\perp}u^n}
  \leq \norm{u^n} \leq \underline{\alpha}^{-1} \Phi(u^n) \leq r_0$,
  the sequence $\seq{P_{U^\perp}u^n}$ is contained in a relative
  compact set in $U^\perp$, hence there is a (strong) accumulation
  point $u^{**} \in U^\perp$. Together with $P_U u^n \rightarrow P_U
  u^*$ we can conclude that there is a subsequence satisfying $u^{n_l}
  \rightarrow P_U u^* + u^{**} = u^{***}$. Moreover, $\seq{u^n}$ is a
  minimizing sequence, so $u^{***}$ has to be a minimizer.
  
  Finally, the whole sequence has to converge to $u^{***}$:
  The mappings $T_n(u) =  J_{s_n}\bigl(u - s_n F'(u) \bigr)$
  satisfy
  \[
  \norm{T_n(u) - T_n(v)} \leq \norm{(I - s_n K^*K)(u - v)} \leq
  \norm{u - v}
  \]
  for all $u,v \in \lpspace{2}$, since all proximal mappings $J_{s_n}$
  are non-expansive and $s_n\leq\tfrac{2}{\norm{K}^2}$.
  So if, for an arbitrary $\varepsilon > 0$ there exists a $n$ such
  that $\norm{u^n - u^{***}} \leq \varepsilon$, then 
  \[
  \norm{u^{n+1} - u^{***}} = \norm{T_n(u^n) - T_n(u^{***})} \leq
  \norm{u^n - u^{***}} \leq \varepsilon
  \]
  since $u^{***}$ is minimizer and hence a fixed point of each $T_n$
  (see Section
  \ref{sec:iter_thres_grad_proj}). By induction, $u^n
  \rightarrow u^{***}$ strongly in $\lpspace{2}$.
\end{Proof}


With the notions of FBI property and strict sparsity pattern from Definitions~\ref{def:fbi} resp.~\ref{def:strict_sparsity_pattern}, one is
able to show linear convergence as soon as one of this two situations
is given.

\begin{Proof}[Proof of
  Theorem~\ref{thm:iter_thres_lin_conv}]
  Observe that the prerequisites of Lemma
  \ref{lem:iter_thres_strong_conv} are fulfilled, so there exists a
  minimizer $u^*$ such that $u^n \rightarrow u^*$ in
  $\lpspace{2}$. Thus, we have to show that each of the two cases
  stated in the theorem leads to a linear convergence rate.
  
  Consider the first case, i.e.~$K$ possesses the FBI property.
  We utilize that, by Lemma~\ref{lem:bregman_inf_dim_est},
  the Bregman-Taylor distance according to
  \eqref{eq:strong_conv_est} and \eqref{eq:taylor_rest} can be estimated
  such that \eqref{eq:rest_local_estimate} is satisfied for some $c > 0$.
  This implies, by Theorem \ref{thm:gen_grad_proj_lin_conv}, the linear
  convergence rate.
  
  For the proof for the second case, 
  we refer to Appendix~\ref{sec:proof_theorem_lin_conv}.
\end{Proof}

  

\begin{Corollary}
  In particular, choosing the step-size constant, i.e.~$s_n = s$ with
  $s \in
  {]0,2\norm{K}^{-2}[}$ also leads to linear convergence under the
  prerequisites of Theorem~\ref{thm:iter_thres_lin_conv},
  for example the step-size $s_n=1$ works for $\norm{K}< \sqrt{2}$.
\end{Corollary}

\begin{Remark}[(Descent and Bregman-Taylor implies linear rate)]
  \label{rem:lin_conv_only_fbi}
  With Theorem~\ref{thm:gen_grad_proj_lin_conv}, the linear
  convergence follows directly from the estimate of the
  Bregman-Taylor distance 
  \[
  \norm{v - u^*}^2 \leq c(M,u^*,K) \bigl( R(v) + T(v) \bigr) \quad
  \text{whenever} \quad (F+\Phi)(v) \leq M
  \]
  which can be established if $K$ satisfies the FBI property. 
  Since the proof of
  Theorem~\ref{thm:gen_grad_proj_lin_conv} relies essentially on
  Proposition~\ref{prop:descent_rate}, one can easily convince oneself
  that the applicability of this proposition is sufficient for linear
  convergence, which is already the case if
  \eqref{eq:rest_distance_est} and $0 < \underline{s} \leq s_n$ 
  is satisfied.
\end{Remark}

\begin{Remark}[(The weak step-size condition as accelerated method)]
  \label{rem:sparsity_weak_step_size}
  As already mentioned in Remark \ref{rem:step_size_descent}, the
  condition on the step-size can be relaxed. In the particular setting
  that $F(u) = \frac12 \norm{Ku - f}^2$, the estimate
  \eqref{eq:weak_step_size_cond} reads as
  \begin{multline*}
    \Bigabs{ \int_0^1 \bigscp{K^*K\bigl(u^n+t(u^{n+1} - u^n)\bigr) - 
        K^*Ku^n}{u^{n+1} - u^n} \dd{t} } \\
    = \frac{\norm{K(u^{n+1} - u^n)}^2}{2} \leq (1 - \delta) D_{s_n}(u^n)
  \end{multline*}
  Now, the choice $s_n$ according to
  \begin{equation}
    \label{eq:sparse_weak_step_size}
    s_n \norm{K(u^{n+1} - u^n)}^2 \leq 2(1-\delta) \norm{u^{n+1} -
      u^n}^2 \ ,
  \end{equation}
  is sufficient for the above, since one has the estimate
  \eqref{eq:dist_norm2}. Together with the boundedness $0 <
  \underline{s} \leq s_n$, this is exactly the step-size `Condition
  (B)' in \cite{daubechies2007projgrad}.
  
  Hence, as can be easily seen, the
  choice gives sufficient descent in order to apply 
  Proposition~\ref{prop:descent_rate}. 
  Consequently, linear convergence
  remains valid for such an `accelerated' iterative soft-thresholding
  procedure if $K$ possesses the FBI property, see
  Remark~\ref{rem:lin_conv_only_fbi}.
\end{Remark}

\begin{Remark}[(Relaxation of the FBI property)]
  %
  It is possible to relax the FBI property.
  Suppose that $K$ fulfills the FBI property \emph{of order $S = \abs{I}$}
  (with the set $I$ defined in the proof of Lemma
  \ref{lem:bregman_inf_dim_est}), 
  i.e., that $K|_I$ is injective for every finite subset $I \subset \NN$ of  
  size less or equal to $S$.
  This
  immediately yields the existence of $\bar c(U,K) > 0$ such
  that
  $\bar c(U,K)\norm{u}^2 \leq \norm{Ku}^2$ for each $u \in U^\perp$ where
  $U^\perp$ is the finite-coefficient subspace as
  defined in the proof of Lemma~\ref{lem:bregman_inf_dim_est}.
  One can easily check that the remaining arguments also remain true
  and consequently,
  Theorem~\ref{thm:iter_thres_lin_conv} still holds.
  
  %
\end{Remark}


The constants in the estimates of Lemma
\ref{lem:bregman_inf_dim_est} and
Theorem~\ref{thm:iter_thres_lin_conv} are in general not
computable unless the solution is determined.
Nonetheless there are some
situations in which prior knowledge about the operator $K$
can be used to estimate the
decay rate.

\begin{Theorem}
  \label{thm:compact_operator}
  Let $K:\lpspace{2}\to\HH_2$, $K \neq 0$ be a compact, linear operator fulfilling the FBI property and define
  \begin{align*}
  \sigma_k &= \inf\ \Bigset{ \frac{\norm{Ku}^2}{\norm{u}^2} }{u \neq 0
    \ , \ u_l = 0 \ \text{for all} \ l \geq k } \ , \\
  \mu_k &= \sup \ \Bigset{ \frac{\norm{Ku}^2}{\norm{u}^2} }{u \neq 0
    \ , \ u_l = 0 \ \text{for all} \ l < k } \ .
  \end{align*}
  Furthermore, choose $k_0$ such that $\mu_{k_0} \leq \underline{\alpha}^2/(4
  \norm{f}^2)$ (with $\infty$ allowed on the right-hand side). 
  Let
  $\{u^n\}$ be a sequence generated by the iterative soft-thresholding algorithm~(\ref{eq:iter_soft_thres}) with initial value $u^0 = 0$ 
  and constant step-size $s= \norm{K}^{-2}$ 
  for the minimization of~(\ref{eq:sparse_min_prob}) and let $u^*$ 
  denote a minimizer.
  Then it holds that
  $\norm{u^n-u^*}\leq C \lambda^n$
  with
  \[
  \lambda = \max\ \Bigl(
  1 - \frac{\sigma_{k_0}}{4\sigma_{k_0}+8\norm{K}^2},
  1 - \frac{\sigma_{k_0}\underline{\alpha}^2}{4\sigma_{k_0}\underline{\alpha}^2 + 2(\sigma_{k_0} + 2\norm{K}^2)\norm{K}^2\norm{f}^2}
  \Bigr)^{1/2}
  \]
  for some $C\geq 0$.
\end{Theorem}
The proof is given in Appendix~\ref{sec:proof-theorem}.



\section{Convergence of related methods}
\label{sec:conv-relat-meth}

In this section, we show how linear convergence can be obtained for
some related methods. In particular,
iterative thresholding methods for minimization problems with
joint sparsity constraints as well as an accelerated gradient 
projection method are considered. 
Both algorithms can be written as a generalized
gradient projection method, hence the analysis carried out in 
Sections~\ref{sec:iter_thres_grad_proj} and \ref{sec:conv_gen_grad_proj} 
can be applied, demonstrating the broad range of applications.

\subsection{Joint sparsity constraints}

First, we consider the situation of so-called joint sparsity for
vector-valued problems, see~\cite{bror2005distributedcompressedsensing,teschke2007nonlinjointsparsity,fornasier2007jointsparsity}.
The problems considered are set in the Hilbert
space $(\lpspace{2})^N$ for some $N \geq 1$ which is interpreted such
that for $u \in (\lpspace{2})^N$ the $k$-th component $u_k$ is a
vector in $\RR^N$. Given a linear and continuous
operator $K: (\lpspace{2})^N \rightarrow
\HH_2$, some data $f \in \HH_2$, a norm $\abs{\placeholder}$ of $\RR^N$
and a sequence $\alpha_k \geq \underline{\alpha} > 0$, the typical
inverse problem with joint sparsity constraints reads as
\begin{equation}
  \label{eq:min_prob_joint_sparse}
  \min_{u \in (\lpspace{2})^N} \ \frac{\norm{Ku - f}^2}{2} +
  \sum_{k=1}^\infty \alpha_k \abs{u_k} \ .
\end{equation}
In many applications, $\abs{\placeholder} = \norm[q]{\placeholder}$
for some $1 \leq q \leq \infty$.

To apply the generalized gradient projection method for
\eqref{eq:min_prob_joint_sparse}, we split the functional into
\[
F(u) = \frac{\norm{Ku-f}^2}{2} \quad , \quad \Phi(u) =
\sum_{k=1}^\infty \alpha_k \abs{u_k} \ .
\]
Analogously to
Proposition~\ref{prop:threshold_gen_grad_proj_equiv}, one needs to
know the associated proximal mappings $J_s$ which can be reduced to
the computation of the proximal
mappings for $\subgrad \abs{\placeholder}$ on $\RR^N$.
These are known to be
\[
(I + s\subgrad\abs{\placeholder})^{-1}(x) = (I -
P_{\sett{\abs{\placeholder}_* \leq s}})(x)
\]
where $P_{\sett{\abs{\placeholder}_* \leq s}}$ denotes the projection to the
closed $s$-ball associated with the dual norm $\abs{\placeholder}_*$.
Again, as can be seen in analogy to
Proposition~\ref{prop:threshold_gen_grad_proj_equiv}, the generalized
gradient projection method for~\eqref{eq:min_prob_joint_sparse} is
given by the iteration
\begin{equation}
  \label{eq:iter_thres_joint_sparse}
  u^{n+1} = \mathbf{S}_{s_n \alpha}\bigl(u^n - s_nK^*(Ku^n - f) \bigr)
  \quad, \quad \bigl(\mathbf{S}_{s_n \alpha}(w)\bigr)_k = (I -
  P_{\sett{\abs{\placeholder}_* \leq s_n\alpha_k}})(w_k) 
\end{equation}
where $\seq{s_n}$ satisfies a suitable step-size rule,
e.g.~according to \eqref{eq:iter_soft_thres_step_size}
or~\eqref{eq:sparse_weak_step_size}.

Let us examine this method with respect to convergence. First, fix a
minimizer $u^*$ which satisfies the optimality condition $w^* =
-K^*(Ku^* - f) \in \subgrad \Phi(u^*)$. As one knows from convex
analysis, this can also be formulated pointwise, and Asplund's
characterization of $\subgrad\abs{\placeholder}$ (see
\cite{showalter1997monotoneoperators}, Proposition II.8.6) leads to
\[
\begin{aligned}
  \abs{w_k^*}_* &\leq \alpha_k  &\text{if}& & u_k^* &= 0 \\
  \abs{w_k^*}_* = \alpha_k \ \text{and} \ &w_k^* \cdot u_k^* = \alpha_k \abs{u_k^*} 
  & \text{if} & & u_k^* & \neq 0
\end{aligned}
\]
where $w_k^* \cdot u_k^*$ denotes the usual inner product of $w_k^*$ and $u_k^*$ 
in $\RR^N$.
Now, one can proceed in complete analogy to the proof of
Lemma~\ref{lem:bregman_inf_dim_est} in order to get an estimate of the
associated Bregman distance: One constructs 
$I = \set{k \in \NN}{\abs{w_k^*}_* = \alpha_k}$ as well as the
closed subspace $U = \set{v \in (\lpspace{2})^N}{v_k = 0 \ \text{if} \
  k \notin I}$
for which $U^\perp$ is finite-dimensional. Furthermore, we have 
$\rho = \sup_{k
  \notin I} \abs{w_k^*}_*/\alpha_k < 1$ and, by equivalence of norms in
$\RR^N$, one gets $C_0, c_0 > 0$ such that $c_0\abs{x}_2 \leq \abs{x}
\leq C_0 \abs{x}_2$ (with $\abs{x}_2^2 = x\cdot x$)
for all $x \in \RR^N$. Then, for a given $M \in \RR$, whenever $(F+\Phi)(v)
\leq M$,
\[
R(v) \geq \frac{(1-\rho) \underline{\alpha}^2
  c_0}{(M+1)C_0} \Bigl( \sum_{k\notin I} \abs{v_k - u_k^*}_2^2
\Bigr)^{1/2} = c_1(M, u^*) \norm{P_U(v - u^*)}^2 \ ,
\]
establishing an analogon of~\eqref{eq:bregman_rest_sparsity_est}.
If $K$ moreover satisfies the FBI property, then one also gets
an analogon to~\eqref{eq:bregman_taylor_rest_sparsity_est}, i.e.
\[
R(v) + T(v) \geq c_2(M,u^*,K) \norm{v - u^*}^2
\]
whenever $(F+\Phi)(v) \leq M$, by arguing
analogously to Lemma~\ref{lem:bregman_inf_dim_est}.

Since these two inequalities are the essential ingredient for proving
convergence as well as the linear rate,
cf.~Lemma~\ref{lem:iter_thres_strong_conv} and
Theorem~\ref{thm:iter_thres_lin_conv}, it holds:

\begin{Theorem}
  The iterative soft-thresholding
  procedure~\eqref{eq:iter_thres_joint_sparse} for the minimization
  problem~\eqref{eq:min_prob_joint_sparse} converges to a minimizer in
  the strong sense in $(\lpspace{2})^N$ if the step-size
  rule~\eqref{eq:iter_soft_thres_step_size} is satisfied.
  
  Furthermore, the convergence will be at linear rate if $K$ possesses
  the FBI property and the step-size rule
  \eqref{eq:sparse_weak_step_size}
  as well as $0 < \underline{s} \leq s_n$ 
  is satisfied. In particular, this is the case when
  $0 < \underline{s} \leq s_n \leq \overline{s} < 2/\norm{K}^2$.
\end{Theorem}

\subsection{Accelerated gradient projection methods}

An alternative approach to implement sparsity constraints for linear
inverse problems is based on minimizing the discrepancy within a
weighted $\lpspace{1}$-ball \cite{daubechies2007projgrad}. With the
notation used in Section
\ref{sec:conv_rates}, the problem can be generally formulated as
\begin{equation}
  \label{eq:acc_grad_proj_func}
  \min_{u \in \Omega} \ \frac{\norm{Ku-f}^2}{2} 
  \quad , \quad \Omega = \Bigset{u \in \lpspace{2}}{\sum_{k=1}^\infty
    \alpha_k \abs{u_k}  \leq 1}  \ . 
\end{equation}
For this classical situation of constrained minimization, one finds
that the generalized gradient projection method and the gradient
projection method coincide (for $F(u)=\tfrac12\norm{Ku-f}^2$ and $\Phi =
I_\Omega$), see Section \ref{sec:iter_thres_grad_proj},
and yield the iteration proposed in
\cite{daubechies2007projgrad}. Consequently, classical convergence
results hold
for a variety of step-size rules \cite{dunn1981gradientprojection},
including the `Condition (B)' introduced in
\cite{daubechies2007projgrad}, see also
Remark~\ref{rem:sparsity_weak_step_size}.

Let us note that linear convergence results can be obtained with the
same techniques which have been used to prove
Theorem~\ref{thm:iter_thres_lin_conv}: First, consider
the Bregman distance $R$ associated with $\Phi = I_\Omega$ in a
minimizer $u^* \in \Omega$. With $w^* = -K^*(Ku^* - f)$,
the optimality condition reads as
\[
\scp{w^*}{v - u^*} \leq 0 \quad \text{for all} \ v \in \Omega \quad
\Leftrightarrow \quad \norm[\infty]{\alpha^{-1}w^*} = \scp{w^*}{u^*} 
\]
where $(\alpha^{-1}w^*)_k = \alpha_k^{-1}w_k^*$ which is in
$\lpspace{\infty}$ since $\alpha_k \geq \underline{\alpha} >
0$. Introduce $I = \set{k}{\abs{\alpha_k^{-1}w_k^*} =
  \norm[\infty]{\alpha^{-1}w^*}}$ which has to be finite since otherwise $w^*
\notin \lpspace{2}$, see the proof of
Lemma~\ref{lem:bregman_inf_dim_est}. Suppose that $w^*\neq 0$ (which
corresponds to $Ku^* \neq f$), so
$\sup_{k \notin I} \
\abs{\alpha_k^{-1}w_k^*}/\norm[\infty]{\alpha^{-1}w^*} = \rho < 1$.
Moreover, 
\begin{equation}
  \label{eq:accel_grad_proj_sol_prop}
  \sum_{k\in I} \alpha_k
  \abs{u_k^*} = 1 \qquad \text{and} \qquad 
  \sign{u_k^*} = \sign{w_k^*} \ \text{for all} \ k \ \text{with} \ u_k^*
  \neq 0 \ , 
\end{equation}
since $\sum_{k\in I} \alpha_k \abs{u_k^*} < 1$ leads to the contradiction
\begin{align*}
  \norm[\infty]{\alpha^{-1} w^*} &=\scp{w^*}{u^*} = \sum_{k=1}^\infty
  \alpha_k^{-1}w_k^* \alpha_k u_k^* 
  \displaybreak[2]
  \\
  &\leq \sum_{k \notin I}
  \abs{\alpha^{-1}_k w_k^*} \abs{\alpha_k u_k^*} + \sum_{k \in I}
  \norm[\infty]{\alpha^{-1}w^*} \alpha_k \abs{u_k^*} \\
  \displaybreak[2]
  &\leq \Bigl( \rho \sum_{k \notin I}\alpha_k \abs{u_k^*} + \sum_{k \in I}
  \alpha_k \abs{u_k^*} \Bigr)
  \norm[\infty]{\alpha^{-1}w^*} < \norm[\infty]{\alpha^{-1}w^*} 
\end{align*}
while $\sign{u_k^*} \neq \sign{w_k^*}$ for some $k$ with $u_k^* \neq 0$
implies the contradiction
\[
\norm[\infty]{\alpha^{-1}w^*} = \sum_{k=1}^\infty \alpha_k^{-1} w_k^*
\alpha_k u_k^* < \sum_{k=1}^\infty \abs{\alpha_k^{-1}w_k^*} \abs{\alpha_k
  u_k^*} \leq \norm[\infty]{\alpha^{-1}w^*} \ .
\]

Moreover, $\sum_{k\in I} \alpha_k \abs{u_k^*} = 1$ also yields $u^*_k
= 0$ for all $k \notin I$. Furthermore, observe that the equation for the
signs in \eqref{eq:accel_grad_proj_sol_prop} gives
$\sum_{k \in I} w_k^*u_k^* = \norm[\infty]{\alpha^{-1}w^*}$.
For $v \notin \Omega$ we have $R(v) = \infty$, so estimate the Bregman
distance for $v \in \Omega$ as follows:
\begin{align*}
  R(v) &= -\scp{w^*}{v - u^*} = \sum_{k\in I} \alpha_k^{-1} w_k^* \alpha_k
  (u_k^* - v_k) - \sum_{k \notin I} \alpha_k^{-1}w_k^* \alpha_k v_k \\
  &\geq
  \norm[\infty]{\alpha^{-1}w^*} - \norm[\infty]{\alpha^{-1}w^*}\sum_{k \in I}
  \alpha_k \abs{v_k} - \rho\norm[\infty]{\alpha^{-1}w^*} \sum_{k \notin I}
  \alpha_k \abs{v_k} \\
  & \geq (1-\rho) \sum_{k\notin I} \alpha_k \abs{v_k} \geq (1-\rho)
  \underline{\alpha}  \norm[1]{P_U v}
\end{align*}
where $U = \set{u \in \lpspace{2}}{u_k = 0 \ \text{for} \ k \in I}$.
Using that $\norm{v} \leq \norm[1]{v}$ as well as 
$\underline{\alpha} \norm{v} \leq \underline{\alpha}
\norm[1]{v} \leq 1$ for all $v \in \Omega$ finally gives, together with
$P_U u^* = 0$,
\[
R(v) \geq (1-\rho) \underline{\alpha}^2 \norm{P_U(v - u^*)}^2 \quad
\text{for all} \ v \in \lpspace{2} \ .
\]
If $K$ possesses the FBI property, one can, analogously to the
argumentation presented in the proof of
Lemma~\ref{lem:bregman_inf_dim_est}, estimate the
Bregman-Taylor distance such that, for some $c(u^*,K) > 0$,
\[
R(v) + T(v) \geq c(u^*,K) \norm{v - u^*}^2 \quad \text{for all} \ v
\in \lpspace{2} \ .
\]
By Theorem~\ref{thm:gen_grad_proj_lin_conv}, the gradient projection
method for \eqref{eq:acc_grad_proj_func} converges linearly. This
remains true for each `accelerated' step-size choice according to
`Condition (B)' in \cite{daubechies2007projgrad}, see
Remark~\ref{rem:sparsity_weak_step_size}.
This result can be summarized in the following theorem.

\begin{Theorem}
  Assume that $K: \lpspace{2} \rightarrow \HH_2$ satisfies the FBI
  property, $\alpha_k \geq \underline{\alpha} > 0$ and $f \in \HH_2
  \setminus K(\Omega)$ where $K(\Omega) = \set{Ku}{\norm[1]{\alpha u} \leq
    1}$. 
  
  Then, the gradient projection method for the minimization problem
  \eqref{eq:acc_grad_proj_func} converges linearly, whenever the
  step-sizes rule~\eqref{eq:sparse_weak_step_size}
  as well as $0 < \underline{s} \leq s_n$ is fulfilled. This is in particular
  the case for  $0 < \underline{s} \leq s_n \leq \overline{s} < 2/\norm{K}^2$.
\end{Theorem}



\section{Conclusions}
\label{sec:conclusions}
We conclude this article with a few remarks on the implications of our results.
We showed that, in many cases, 
iterative soft-thresholding algorithms converge
with linear rate and moreover that there are situations in which
the constants can be calculated explicitly, see 
Theorem~\ref{thm:compact_operator}.
In general, however, the factor $\lambda$, 
which determines the speed within the class of
linearly-convergent algorithms, always depends on the operator $K$ 
but in the considered cases also on the 
initial value $u^0$ and a solution $u^*$. Unfortunately, the dependence
on a solution can cause $\lambda$ to be arbitrarily close to $1$,
meaning that the iterative soft-thresholding converges arbitrarily
slow in some sense, which is also often observed in practice. 

One key ingredient for proving the convergence result 
is the FBI property. This property
also plays a role in the performance analysis of Newton methods 
applied to minimization problems with sparsity 
constraints~\cite{griesse2007ssnsparsity}
and error estimates for $\ell^1$-regularization~\cite{lorenz2008reglp}.
As we have moreover seen, linear convergence can
also be obtained whenever we have convergence a solution with 
strict sparsity pattern. This result is closely connected 
with the fact that~\eqref{eq:sparse_min_prob}, considered on a fixed sign 
pattern, is a quadratic problem, and hence the iteration becomes linear from
some index on. The latter observation is also basis of a couple of 
different algorithms \cite{osborne2000variableselection,efron2004lars,figueiredo2007gradproj}.

At last we want to remark that Theorem~\ref{thm:gen_grad_proj_lin_conv}
on linear convergence of the generalized gradient projection method
holds in general and has been applied in a special case in order to 
prove Theorem~\ref{thm:iter_thres_lin_conv}.
This generality also allowed for a unified treatment of the similar algorithms
presented in Section~\ref{sec:conv-relat-meth} as well as other penalty 
terms such as powers of certain $2$-convex norms, see 
Remark~\ref{rem:bregman_sufficient}. 
In all of these situations, linear convergence
follows from descent properties on the one hand and Bregman (-Taylor) 
estimates on the other hand.

\appendix

\section{Proof of Theorem~\ref{thm:iter_thres_lin_conv} (continued)}
\label{sec:proof_theorem_lin_conv}

For the second case, let $u^*$ possess a strict sparsity pattern.
Define, analogously to the above, the subspace 
$U = \set{v \in \lpspace{2}}{v_k =
  0 \ \text{if} \ u_k^* \neq 0}$. The desired 
result then is implied by the fact that there is an $n_0$ such that
each $u^{n+1}$ with $n \geq n_0$ can be written as
\[
u^{n+1} = (I - s_nP_{U^\perp}K^*KP_{U^\perp})(u^n - u^*) + u^* \ .
\]
For this purpose, we introduce the notations $w^n = -K^*(Ku^n - f)$,
$w^* = -K^*(Ku^* - f)$ and recall the optimality condition $w^*
\in \subgrad \Phi(u^*)$ which can be written as
\begin{align*}
  w_k^* &\in [-\alpha_k,\alpha_k] & \text{if} \ u_k^* &= 0 
  \displaybreak[2]
  \\
  w_k^* &= \alpha_k & \text{if} \ u_k^* &> 0 
  \displaybreak[2]
  \\
  w_k^* &= -\alpha_k & \text{if} \ u_k^* &< 0 \ .
\end{align*}
Due to assumption that $u^*$ has a strict sparsity pattern, $w_k^*
\in {]{-\alpha_k, \alpha_k}[}$ if $u_k^* = 0$,
and hence there is a $\rho > 0$ such that
\begin{align*}
  w_k^* &\in [-(1-\rho)\alpha_k,(1-\rho)\alpha_k] &\text{if} \ u_k^*
  &= 0 
\end{align*}
since $w_k^* \rightarrow 0$ for $k \rightarrow \infty$.
Also note that $u^n \rightarrow u^*$ implies $w^n \rightarrow w^*$
and especially pointwise convergence.

We will treat each of the cases $u_k^* = 0$, $u_k^* > 0$ and $u_k^*
< 0$ separately.

\noindent\textbf{The case $u^*_k=0$:}
First, we find an index $n_1$ 
such that, for $n \geq n_1$,
\[
\norm{u^n - u^*} \leq
\frac{\rho}2\underline{\alpha}\,\underline{s} \quad , \quad
\norm{w^n - w^*} \leq \frac{\rho}2\underline{\alpha}.
\]
So, if $k \in I_0$ with $I_0 = \set{k}{u_k^* = 0}$, we have
\[
\abs{u_k^n} \leq \frac{\rho}{2} s_n \alpha_k \quad , \quad
\abs{w_k^n} \leq \abs{w_k^*} + \abs{w_k^n - w_k^*} \leq
(1-\rho)\alpha_k + \frac{\rho}2\alpha_k
\]
for each $n \geq n_1$. Consequently, for all of these $k$ and $n$,
\[
\abs{u^n - s_n K^*(Ku^n - f)}_k \leq s_n\alpha_k
\]
hence the thresholding operation according to
\eqref{eq:iter_soft_thres} gives $u^{n+1}_k = 0$ for all $n \geq
n_1$ and all $k \in I_0$. Thus, the iteration for $P_Uu^n$ can be
expressed by
\begin{equation}
  \label{eq:proof_lin_conv_zero_set}
  P_U u^{n+1} = 
  (I - s_nP_{U^\perp}K^*KP_{U^\perp})(u^n - u^*) + P_Uu^*
\end{equation}
for all $n \geq n_1 + 1$ since $P_Uu^n = P_U u^* = 0$.

\noindent\textbf{The case $u^*_k>0$:}
Next, investigate all $k \in I_+$ with $I_+ = \set{k}{u_k^* > 0}$. 
This has to be a finite set, so there is a $\delta_+ \in
{]{0,\underline{\alpha}}[}$
such that $u_k^* \geq \delta_+$ for each of such $k$. 
So, choose $n_+$
according to the requirements that for all $n \geq n_+$
\[
\norm{u^n - u^*} \leq \frac{\delta_+}2 \quad , \quad \norm{w^n -
  w^*} \leq \frac{\delta_+}{2\overline{s}} \ .
\]
Then, remembering that $w_k^* = \alpha_k$,
\begin{align*}
  u^n_k + s_nw^n_k &= u^*_k + u^n_k - u^*_k + s_n(w^n_k - w^*_k) +
  s_nw^*_k \\
  &\geq u^*_k - \abs{u^n_k - u^*_k} - s_n\abs{w^n_k - w^*_k} + s_n
  \alpha_k \\
  & \geq \delta_+ - \frac{\delta_+}2 - \frac{s_n
    \delta_+}{2\overline{s}} + s_n\alpha_k 
  \geq s_n \alpha_k
\end{align*}
and hence the iteration gives, by $(w^n - w^*) = -K^*K(u^n - u^*)$,
\begin{align}
  \label{eq:proof_lin_conv_plus_set}
  \notag
  u^{n+1}_k &= u^n_k + s_n w^n_k - s_n\alpha_k \\
  \notag
  &= u^n_k - u^*_k + s_n(w^n_k - w^*_k) + u^*_k \\
  &= \bigl((I - s_nK^*K)(u^n - u^*)\bigr)_k
  + u^*_k
\end{align}
for all $n \geq n_+$ and all $k \in I_+$.

\noindent\textbf{The case $u^*_k<0$:}
Analogously, considering the indices $k \in I_-$ with $I_- =
\set{k}{u_k^* < 0}$,
one can find an $n_-$ such that 
\begin{equation}
  \label{eq:proof_lin_conv_minus_set}
  u^{n+1}_k = \bigl((I - s_nK^*K)(u^n - u^*)\bigr)_k + u^*_k
\end{equation}
also holds for all $n \geq n_-$ and all $k \in I_-$.

Choosing $n_0 = \max\ (n_1+1, n_+, n_-)$ and considering
\eqref{eq:proof_lin_conv_zero_set}--\eqref{eq:proof_lin_conv_minus_set}
as well as remembering that $P_Uu^n = 0$ for $n>n_0$ yields that
indeed
\begin{equation}
  \label{eq:proof_lin_conv_lin_iter}
  u^{n+1} = (I - s_n P_{U^\perp}K^*KP_{U^\perp})(u^n - u^*) + u^* \ .
\end{equation}

Eventually, we can split the iteration into the subspaces $V =
\kernel{KP_{U^\perp}}$ and $V^\perp$, where $V^\perp$ is taken with respect 
to $U^\perp$. For $n \geq n_0$, 
\[
P_V u^{n+1} = (P_V - s_n P_VP_{U^\perp}K^*KP_{U^\perp})(u^n - u^*) +
P_Vu^* = P_Vu^n
\]
due to the fact that $V = \kernel{KP_{U^\perp}} =
\range{P_{U^\perp}K^*}^{\perp}$. Consequently, $P_Vu^n = P_Vu^*$
since there would not hold that $u^n \rightarrow u^*$ otherwise.
Note that $V^\perp$ is finite dimensional, hence there is a $c > 0$
such that $c\norm{P_{V^\perp}u}^2 \leq \norm{KP_{U^\perp}P_{V^\perp}u}^2 =
\norm{KP_{V^\perp}u}^2$ for all $u \in \lpspace{2}$. 
Consequently, each of the self-adjoint mappings
$P_{V^\perp} - s_nP_{V^\perp}K^*KP_{V^\perp}$ is a strict
contraction on $V^\perp$:
\begin{align*}
  \sup_{\norm{P_{V^\perp}u} = 1} \ \bigl| \langle(P_{V^\perp} &-
  s_nP_{V^\perp}K^*KP_{V^\perp})u, \ P_{V^\perp}u \rangle \bigr| \\ 
  &= \sup_{\norm{P_{V^\perp}u} = 1} \ \bigabs{\norm{P_{V^\perp}u}^2 -
    s_n \norm{KP_{V^\perp}}^2} \\ 
  &\leq
  \max \ \Bigl( \overline{s} \norm{K}^2 - 1,
  \sup_{\norm{P_{V^\perp}u} = 1} \ \norm{P_{V^\perp}u}^2 - s_n
  c\norm{P_{V^\perp}u}^2 \Bigr) \\
  &\leq \max \ \bigl( \overline{s}\norm{K}^2 - 1, 1 - \underline{s}c
  \bigr) = \lambda < 1 \ .
\end{align*}
Using that
$u^n - u^* = P_{V^\perp}(u^n - u^*)$ for $n \geq
n_0$ gives, plugged into \eqref{eq:proof_lin_conv_lin_iter},
\[
u^{n+1} - u^* = (P_{V^\perp} -
s_nP_{V^\perp}K^*KP_{V^\perp})(u^n - u^*)
\]
so
\[
\norm{u^{n+1} - u^*}^2 = \norm{P_{V^\perp}(u^{n+1} - u^*)}^2  \leq
\lambda^2 \norm{P_{V^\perp}(u^{n} - u^*)}^2 = \lambda^2  \norm{u^n -
  u^*}^2 \ ,
\]
meaning $\norm{u^n - u^*} \leq \lambda^{n-n_0} \norm{u^{n_0} -
  u^*}$ for $n \geq n_0$. 
Finally, it is easy to find a $C > 0$ such that $\norm{u^n
  - u^*} \leq C\lambda^n$ for all $n$.

\section{Proof of Theorem~\ref{thm:compact_operator}}
\label{sec:proof-theorem}
 \begin{Proof}[Proof of Theorem~\ref{thm:compact_operator}]
  Note that $\sigma_k > 0$ because of the FBI property and that
  $\mu_k \rightarrow 0$ as $k \rightarrow
  \infty$ since $K$ is compact
  (otherwise there would be a bounded sequence which converges
  weakly to zero with images not converging in the strong sense).

  Our aim is to compute a constant $c_1 > 0$ such that $c_1 \norm{P_k(v -
    u^*)}^2 \leq R(v)$ on a suitable bounded set and for a suitable
  $k$. Here, 
  $P_k$ denotes the orthogonal projection onto the subspace
  $\set{u \in \lpspace{2}}{u_l = 0 \ \text{for}\ l < k}$.
  We can assume without loss of generality that $f \neq 0$ and 
  thus estimate the norm of $Ku^* - f$:
  \[
  \frac{\norm{Ku^* - f}^2}{2} \leq (F+\Phi)(u^*) \leq (F+\Phi)(0) =
  \frac{\norm{f}^2}{2} \quad \Rightarrow \quad \norm{Ku^* - f} \leq
  \norm{f}
  \]
  Because the index $k_0$ is chosen such that $\mu_{k_0} \leq \underline{\alpha}^2/(4 \norm{f}^2)$ we can estimate
  \begin{align*}
    \norm{P_{k_0}K^*(Ku^* - f)} &= \sup_{\norm{v} \leq 1} \ \scp{Ku^* -
      f}{KP_{k_0} v} \\
    &\leq \sup_{\norm{v} \leq 1} \ \norm{Ku^* - f}
    \norm{KP_{k_0} v} \leq \mu_{k_0}^{1/2} \norm{f} \leq
    \underline{\alpha}/2
  \end{align*}
  and consequently, $w^* = -K^*(Ku^* - f)$ satisfies $\abs{w_k^*} \leq
  \underline{\alpha}/2$ for each $k \geq k_0$. 
  Recall from the proof of Lemma~\ref{lem:bregman_inf_dim_est} that
  this in particular means that $u^*_k = 0$, so one obtains the
  estimate
  \[
  R(v) \geq \sum_{k \geq k_0} \alpha_k(\abs{v_k} - \abs{u_k^*}) +
  w_k^* v_k
  \geq \tfrac12 \underline{\alpha} \sum_{k \geq k_0} \abs{v_k - u_k^*}
  \geq \tfrac12 \underline{\alpha} \norm{P_{k_0}(v - u^*)} \ .
  \]
  We assumed that the first iterate is $u^0 = 0$, 
  so $(F+\Phi)(u^n) \leq \frac{\norm{f}^2}{2}$. Consequently,
  $\norm{P_{k_0}v}^{-1} \geq 2\underline{\alpha}\norm{f}^{-2}$ whenever
  $\Phi(v) \leq (F+\Phi)(0)$, so
  \[
  R(v) \geq \underline{\alpha}^2\norm{f}^{-2} \norm{P_{k_0}(v - u^*)}^2 \ .
  \]
  An estimate for the Taylor-distance $T$ is found with the
  help of $\sigma_{k_0}$:
  \begin{align*}
    T(v) &= \frac{\norm{K(v - u^*)}^2}{2} \geq
    \frac{\norm{KP_{k_0}^\perp(v - u^*)}^2}{4} - \frac{\norm{K}^2
      \norm{P_{k_0}(v - u^*)}^2}{2} \displaybreak[2]\\
    & \geq \frac{\sigma_{k_0}}{4} \bigl( \norm{P_{k_0}^\perp(v - u^*)}^2
    + \norm{P_{k_0}(v-u^*)}^2 \bigr) 
    -  \frac{(\sigma_{k_0} + 2 \norm{K}^2)
      \norm{f}^2}{4\underline{\alpha}^2} R(v) \ ,
  \end{align*}
  where $P_{k_0}^\perp = I - P_{k_0}$.
  Rearranging terms gives:
  \[
  \frac{\sigma_{k_0}}{4} \norm{v - u^*}^2 \leq \max \ \Bigl (1,
  \frac{(\sigma_{k_0} + 2\norm{K}^2)\norm{f}^2}{4\underline{\alpha}^2}
  \Bigr) \bigl( R(v) + T(v) \bigr)
  \]
  leading to the desired constant $c$ in Proposition
  \ref{prop:descent_rate}:
  \[
  \norm{v - u^*}^2 \leq \max \ \Bigl( \frac{4}{\sigma_{k_0}}, 
  \frac{(\sigma_{k_0} + 2\norm{K}^2)\norm{f}^2}
  {\sigma_{k_0} \underline{\alpha}^2} \Bigr)
  r_n \ ,
  \]
  namely $c = \max\ \Bigl( \frac{4}{\sigma_{k_0}}, 
  \frac{(\sigma_{k_0} + 2\norm{K}^2)\norm{f}^2}
  {\sigma_{k_0} \underline{\alpha}^2} \Bigr)$.
  Estimating $\lambda$ according to \eqref{eq:rest_exp_decay}
  with constant step-size $s = \norm{K}^{-2}$ and 
  $\delta = 1 - s \norm{K}^2/2 = 1/2$ 
  yields 
  \begin{align*}
  \lambda^2 & \leq 1 - \frac{1}{4  + 2c\norm{K}^2} \\
  & = \max \ \Bigl(
    1 - \frac{\sigma_{k_0}}{4\sigma_{k_0}+8\norm{K}^2},
    1 - \frac{\sigma_{k_0}\underline{\alpha}^2}{4\sigma_{k_0}\underline{\alpha}^2 + 2(\sigma_{k_0} + 2\norm{K}^2)\norm{K}^2\norm{f}^2}
    \Bigr) \ .
  \end{align*}
 \end{Proof}
 \begin{Remark}
  The proof of Proposition \ref{prop:descent_rate} also establishes
  $\norm{u^n - u^*} \leq (cr_0)^{1/2} \lambda^n$ which
  implies in turn, by estimating $r_0 \leq (F+\Phi)(0) = \norm{f}^2/2$
  and the maximum by the sum, the a-priori estimate 
  \begin{multline*}
    \label{eq:a_priori_norm_est}
    \norm{u^n - u^*} \leq 
    \sqrt{ 
      \frac{4 \underline{\alpha}^2\norm{f}^2 + (\sigma_{k_0} + 2
        \norm{K}^2)\norm{f}^4 }{2 \sigma_{k_0}
        \underline{\alpha}^2} } \\
    \cdot
    \max \ \Bigl(
    1 - \frac{\sigma_{k_0}}{4\sigma_{k_0}+8\norm{K}^2},
    1 - \frac{\sigma_{k_0}\underline{\alpha}^2}{4\sigma_{k_0}\underline{\alpha}^2 + 2(\sigma_{k_0} + 2\norm{K}^2)\norm{K}^2\norm{f}^2}
    \Bigr)^{n/2} \ .
  \end{multline*}
 \end{Remark}

\section*{References}

\bibliography{literature,paper}
\bibliographystyle{plain}

{\footnotesize
\centerline{\rule{9pc}{.01in}}
\bigskip
\centerline{
  Center for Industrial Mathematics / Fachbereich 3}
\centerline{
  University of Bremen, Postfach 33 04 40}
\centerline{
  28334 Bremen, Germany}
\centerline{e-mail: kbredies@math.uni-bremen.de, dlorenz@math.uni-bremen.de}
}

\end{document}